# A CONVERSE OF THE JENSEN INEQUALITY FOR CONVEX MAPPINGS OF SEVERAL VARIABLES AND APPLICATIONS

S. S. DRAGOMIR

ABSTRACT. In this paper we point out a converse result of the celebrated Jensen inequality for differentiable convex mappings of several variables and apply it to counterpart well-known analytic inequalities. Applications to Shannon's and Rényi's entropy mappings are also given.

## 1. INTRODUCTION

Let $f : X \to \mathbb{R}$ be a convex mapping defined on the linear space $X$ and $x_i \in X$, $p_i \geq 0$ $(i = 1, ..., m)$ with $P_m := \sum_{i=1}^{m} p_i > 0$.

The following inequality is well known in the literature as Jensen's inequality

$$(1.1) \qquad f\left(\frac{1}{P_m} \sum_{i=1}^{m} p_i x_i\right) \leq \frac{1}{P_m} \sum_{i=1}^{m} p_i f(x_i).$$

There are many well known inequalities which are particular cases of Jensen's inequality such as the weighted arithmetic mean-geometric mean- harmonic mean inequality, the Ky Fan inequality, the Hölder inequality, etc. For a comprehensive list of recent results on the Jensen inequality, see the book [1] and the papers [2]-[14] where further results are given

In this paper, we point out a converse inequality for Jensen's result (1.1) for the case of differentiable convex mappings whose partial derivatives are bounded. Applications for some particular inequalities and for the Shannon and Rényi entropy mappings are also given.

## 2. A CONVERSE INEQUALITY

The following converse of Jensen's inequality holds.

**Theorem 1.** *Let $f : \mathbb{R}^n \to \mathbb{R}$ be a differentiable convex mapping and $x_i \in \mathbb{R}^n$, $i = 1, ..., m$. Suppose that there exist the vectors $\psi, \phi \in \mathbb{R}^n$ such that $\psi \leq x_i \leq \phi$ (the order is considered on the co-ordinates) and $m, M \in \mathbb{R}^n$ are such that $m \leq \frac{\partial f(x_i)}{\partial x} \leq M$ for all $i \in \{1, ..., m\}$. Then for all $p_i \geq 0$ $(i = 1, ..., m)$ with $P_m := \sum_{i=1}^{m} p_i > 0$, we have the inequality*

$$(2.1) \qquad 0 \leq \frac{1}{P_m} \sum_{i=1}^{m} p_i f(x_i) - f\left(\frac{1}{P_m} \sum_{i=1}^{m} p_i x_i\right) \leq \frac{1}{4} \|\Phi - \phi\| \|M - m\|,$$

*where $\|\cdot\|$ is the usual Euclidean norm on $\mathbb{R}^n$.*

*Proof.* Firstly, and for the sake of completeness, we prove the following inequality for convex functions obtained by Dragomir and Goh in [14]:

$$
\begin{align}
(2.2) \quad 0 &\leq \frac{1}{P_m} \sum_{i=1}^{m} p_i f(x_i) - f\left(\frac{1}{P_m} \sum_{i=1}^{m} p_i x_i\right) \\
&\leq \frac{1}{P_m} \sum_{i=1}^{m} p_i \langle x_i, \nabla f(x_i) \rangle - \left\langle \frac{1}{P_m} \sum_{i=1}^{m} p_i x_i, \frac{1}{P_m} \sum_{i=1}^{m} p_i \nabla f(x_i) \right\rangle,
\end{align}
$$

where $\langle \cdot, \cdot \rangle$ is the usual inner product on $\mathbb{R}^n$ and

$$\nabla f(x) = \frac{\partial f(x)}{\partial x} := \left(\frac{\partial f(x)}{\partial x^1}, ..., \frac{\partial f(x)}{\partial x^n}\right)$$

is the vector of the partial derivatives, $x = (x^1, ..., x^n) \in \mathbb{R}^n$.

As $f : \mathbb{R}^n \to \mathbb{R}$ is differentiable convex, we have the inequality

$$(2.3) \qquad f(x) - f(y) \geq \langle \nabla f(y), x - y \rangle, \text{ for all } x, y \in \mathbb{R}^n.$$

Choose in (2.3), $x = \frac{1}{P_m} \sum_{i=1}^{m} p_i x_i$ and $y = x_j$ to obtain

$$(2.4) \qquad f\left(\frac{1}{P_m} \sum_{i=1}^{m} p_i x_i\right) - f(x_j) \geq \left\langle \nabla f(x_j), \frac{1}{P_m} \sum_{i=1}^{m} p_i x_i - x_j \right\rangle$$

for all $j \in \{1, ..., n\}$.

If we multiply (2.4) by $p_j \geq 0$ and sum over $j$ from 1 to $m$, we obtain

$$
\begin{align}
P_m f&\left(\frac{1}{P_m} \sum_{i=1}^{m} p_i x_i\right) - \sum_{j=1}^{m} p_j f(x_j) \\
&\geq \frac{1}{P_m} \left\langle \sum_{j=1}^{m} p_j \nabla f(x_j), \sum_{i=1}^{m} p_i x_i \right\rangle - \sum_{j=1}^{m} \langle \nabla f(x_j), x_j \rangle.
\end{align}
$$

Dividing by $P_m > 0$, we obtain (2.2).

A simple calculation shows that

$$
\begin{align}
(2.5) \quad &\frac{1}{P_m} \sum_{i=1}^{m} p_i \langle x_i, \nabla f(x_i) \rangle - \left\langle \frac{1}{P_m} \sum_{i=1}^{m} p_i x_i, \frac{1}{P_m} \sum_{i=1}^{m} p_i \nabla f(x_i) \right\rangle \\
&= \frac{1}{2P_m^2} \sum_{i,j=1}^{m} p_i p_j \langle x_i - x_j, \nabla f(x_i) - \nabla f(x_j) \rangle.
\end{align}
$$

Taking the modulus in both parts of (2.5), and accounting for the fact that the left hand side is positive (by (2.2)), we obtain, by Schwartz's inequality in inner product spaces, i.e., we recall it $|\langle a, b \rangle| \leq \|a\| \|b\|$, $a, b \in \mathbb{R}^n$, that

$$
\begin{align}
(2.6) \quad &\frac{1}{P_m} \sum_{i=1}^{m} p_i \langle x_i, \nabla f(x_i) \rangle - \left\langle \frac{1}{P_m} \sum_{i=1}^{m} p_i x_i, \frac{1}{P_m} \sum_{i=1}^{m} p_i \nabla f(x_i) \right\rangle \\
&\leq \frac{1}{2P_m^2} \sum_{i,j=1}^{m} p_i p_j |\langle x_i - x_j, \nabla f(x_i) - \nabla f(x_j) \rangle| \\
&\leq \frac{1}{2P_m^2} \sum_{i,j=1}^{m} p_i p_j \|x_i - x_j\| \|\nabla f(x_i) - \nabla f(x_j)\|.
\end{align}
$$

Using the Cauchy-Buniakowsky-Schwartz inequality for double sums, we can state that

$$(2.7) \quad \frac{1}{2P_m^2} \sum_{i,j=1}^{m} p_i p_j \|x_i - x_j\| \|\nabla f(x_i) - \nabla f(x_j)\|$$

$$\leq \left( \frac{1}{2P_m^2} \sum_{i,j=1}^{m} p_i p_j \|x_i - x_j\|^2 \right)^{\frac{1}{2}} \times \left( \frac{1}{2P_m^2} \sum_{i,j=1}^{m} p_i p_j \|\nabla f(x_i) - \nabla f(x_j)\|^2 \right)^{\frac{1}{2}}.$$

As a simple calculation shows that

$$\frac{1}{2P_m^2} \sum_{i,j=1}^{m} p_i p_j \|x_i - x_j\|^2 = \frac{1}{P_m} \sum_{i=1}^{m} p_i \|x_i\|^2 - \left\| \frac{1}{P_m} \sum_{i=1}^{m} p_i x_i \right\|^2$$

and

$$\frac{1}{2P_m^2} \sum_{i,j=1}^{m} p_i p_j \|\nabla f(x_i) - \nabla f(x_j)\|^2$$

$$= \frac{1}{P_m} \sum_{i=1}^{m} p_i \|\nabla f(x_i)\|^2 - \left\| \frac{1}{P_m} \sum_{i=1}^{m} \nabla f(x_i) \right\|^2,$$

we can state, by (2.6) and (2.7), that

$$(2.8) \quad \frac{1}{P_m} \sum_{i=1}^{m} p_i \langle x_i, \nabla f(x_i) \rangle - \left\langle \frac{1}{P_m} \sum_{i=1}^{m} p_i x_i, \frac{1}{P_m} \sum_{i=1}^{m} p_i \nabla f(x_i) \right\rangle$$

$$\leq \left( \frac{1}{P_m} \sum_{i=1}^{m} p_i \|x_i\|^2 - \left\| \frac{1}{P_m} \sum_{i=1}^{m} p_i x_i \right\|^2 \right)^{\frac{1}{2}}$$

$$\times \left( \frac{1}{P_m} \sum_{i=1}^{m} p_i \|\nabla f(x_i)\|^2 - \left\| \frac{1}{P_m} \sum_{i=1}^{m} \nabla f(x_i) \right\|^2 \right)^{\frac{1}{2}}.$$

Now, let us observe that, by a simple calculation

$$(2.9) \quad \frac{1}{P_m} \sum_{i=1}^{m} p_i \|x_i\|^2 - \left\| \frac{1}{P_m} \sum_{i=1}^{m} p_i x_i \right\|^2$$

$$= \left\langle \phi - \frac{1}{P_m} \sum_{i=1}^{m} p_i x_i, \frac{1}{P_m} \sum_{i=1}^{m} p_i x_i - \psi \right\rangle - \frac{1}{P_m} \sum_{i=1}^{m} p_i \langle \phi - x_i, x_i - \psi \rangle.$$

As $\psi \leq x_i \leq \phi$ $(i \in \{1, ..., m\})$, then $\langle \phi - x_i, x_i - \psi \rangle \geq 0$ for all $i \in \{1, ..., m\}$ and then

$$\sum_{i=1}^{m} p_i \langle \phi - x_i, x_i - \psi \rangle \geq 0$$

and, by (2.9), we obtain

$$
(2.10) \quad \frac{1}{P_m} \sum_{i=1}^m p_i \|x_i\|^2 - \left\| \frac{1}{P_m} \sum_{i=1}^m p_i x_i \right\|^2
$$
$$
\leq \left\langle \phi - \frac{1}{P_m} \sum_{i=1}^m p_i x_i, \frac{1}{P_m} \sum_{i=1}^m p_i x_i - \psi \right\rangle.
$$

It is known that if $y, z \in \mathbb{R}^n$, then

$$
(2.11) \quad 4 \langle z, y \rangle \leq \|z + y\|^2
$$

with equality iff $z = y$.

Now, if we apply (2.11) for the vectors $z = \phi - \frac{1}{P_M} \sum_{i=1}^m p_i x_i$, $y = \frac{1}{P_m} \sum_{i=1}^m p_i x_i - \psi$, we deduce

$$
\left\langle \phi - \frac{1}{P_m} \sum_{i=1}^m p_i x_i, \frac{1}{P_m} \sum_{i=1}^m p_i x_i - \psi \right\rangle \leq \frac{1}{4} \|\phi - \psi\|^2
$$

and then, by (2.9) - (2.10), we deduce that

$$
(2.12) \quad \frac{1}{P_m} \sum_{i=1}^m p_i \|x_i\|^2 - \left\| \frac{1}{P_m} \sum_{i=1}^m p_i x_i \right\|^2 \leq \frac{1}{4} \|\phi - \psi\|^2.
$$

Similarly, we can state that

$$
(2.13) \quad \frac{1}{P_m} \sum_{i=1}^m p_i \|\nabla f(x_i)\|^2 - \left\| \frac{1}{P_m} \sum_{i=1}^m \nabla f(x_i) \right\|^2 \leq \frac{1}{4} \|M - m\|^2.
$$

Finally, by (2.8) and (2.12) - (2.13), we deduce

$$
(2.14) \quad \frac{1}{P_m} \sum_{i=1}^m p_i \langle x_i, \nabla f(x_i) \rangle - \left\langle \frac{1}{P_m} \sum_{i=1}^m p_i x_i, \frac{1}{P_m} \sum_{i=1}^m p_i \nabla f(x_i) \right\rangle
$$
$$
\leq \frac{1}{4} \|\Phi - \phi\| \|M - m\|,
$$

which, by (2.2), gives the desired inequality (2.1). ∎

**Remark 1.** *A similar result for integrals can be stated, but we omit the details.*

**Remark 2.** *The conditions*

$$
(2.15) \quad \psi \leq x_i \leq \phi, \ m \leq \frac{\partial f(x_i)}{\partial x} \leq M \ (i = 1, ..., m)
$$

*can be replaced by the more general conditions*

$$
(2.16) \quad \sum_{i=1}^m p_i \langle \phi - x_i, x_i - \psi \rangle \geq 0 \ \text{and} \ \sum_{i=1}^m p_i \left\langle M - \frac{\partial f(x_i)}{\partial x}, \frac{\partial f(x_i)}{\partial x} - m \right\rangle \geq 0
$$

*and the conclusion (2.1) will still be valid.*

**Remark 3.** *Even if the new inequality (2.1) is not as sharp as the Dragomir-Goh inequality (2.2), it may be more useful in practice when only some bounds of the* partial derivatives $\frac{\partial f}{\partial x}$ and of the vectors $x_i$ $(i = 1, ..., m)$ are known. Namely, it provides the opportunity to estimate the difference

$$\frac{1}{P_m} \sum_{i=1}^{m} p_i f(x_i) - f\left(\frac{1}{P_m} \sum_{i=1}^{m} p_i x_i\right) =: \Delta(f, x, p)$$

when the quantities $\|\phi - \psi\|$ and $\|M - m\|$ are known.

For example, if the partial derivatives $\frac{\partial f}{\partial x}$ are bounded, i.e., there exists $m, M \in \mathbb{R}^n$ such that $m \leq \frac{\partial f}{\partial x} \leq M$ on the co-ordinates, and if we choose the vector $x_i$ $(i = 1, ..., m)$ not "very far" from a constant vector $x_0$, i.e., $\|\phi - \psi\| \leq \frac{4\varepsilon}{\|M-m\|}$, $\varepsilon > 0$, then, by (2.1), we can conclude that

$$0 \leq \Delta(f, x, p) \leq \varepsilon.$$

The case of convex mappings of a real variable can be stated as follows [20].

**Corollary 1.** *Let $f : \mathbb{R} \to \mathbb{R}$ be a differentiable convex mapping and $x_i \in I$ for all $i \in \{1, ..., m\}$. Then we have the inequality:*

$$\begin{aligned} (2.17) \qquad 0 &\leq \frac{1}{P_m} \sum_{i=1}^{m} p_i f(x_i) - f\left(\frac{1}{P_m} \sum_{i=1}^{m} p_i x_i\right) \\ &\leq \frac{1}{4}(M - m)(f'(M) - f'(m)), \end{aligned}$$

*where $p_i > 0$ $(i = 1, ..., m)$ and $P_m := \sum_{i=1}^{m} p_i > 0$.*

The proof is obvious by the above findings, taking into account that the mapping $f'$ is monotonic nondecreasing, and then $f'(m) \leq f'(x_i) \leq f'(M)$ for all $i \in \{1, ..., m\}$.

3. Applications for Weighted Means

Consider the classical weighted means:

$$\begin{aligned} A_n(\bar{p}, x) &:= \sum_{i=1}^{n} p_i x_i \quad \text{-} \quad \text{the arithmetic mean,} \\ G_n(\bar{p}, x) &:= \prod_{i=1}^{n} x_i^{p_i} \quad \text{-} \quad \text{the geometric mean,} \\ H_n(\bar{p}, x) &:= \frac{1}{\sum_{i=1}^{n} \frac{p_i}{x_i}} \quad \text{-} \quad \text{the harmonic mean,} \end{aligned}$$

provided that $x_i > 0$ $\left(i = \overline{1, n}\right)$ and $p_i$ $(i = 1, ..., n)$ is a probability distribution, i.e., $\sum_{i=1}^{n} p_i = 1$.

The following inequality is well-known in the literature as the *arithmetic mean-geometric mean-harmonic mean inequality*

$$(3.1) \qquad A_n(\bar{p}, x) \geq G_n(\bar{p}, x) \geq H_n(\bar{p}, x),$$

with equality iff $x_1 = ... = x_n$ (for $p_i > 0$, $i = 1, ..., n$).

The following corollaries hold.

**Corollary 2.** *Let $0 < m \leq x_i \leq M < \infty$, $p_i > 0$ $(i = 1, ..., n)$. Then*

$$(3.2) \qquad 1 \leq \frac{A_n(\bar{p}, x)}{G_n(\bar{p}, x)} \leq \exp\left[\frac{(M-m)^2}{4mM}\right].$$

*Equality holds in (3.2) simultaneously iff $x_1 = ... = x_n$.*

The proof follows by the inequality (2.17) for $f(x) = -\ln x$, $x > 0$.

**Corollary 3.** *Let $0 < m \leq y_i \leq M < \infty$, $p_i > 0$ $(i = 1, ..., n)$. Then*

$$(3.3) \qquad 1 \leq \frac{G_n(\bar{p}, \bar{y})}{H_n(\bar{p}, \bar{y})} \leq \exp\left[\frac{(M-m)^2}{4mM}\right].$$

*Equality holds iff $y_1 = ... = y_n$.*

**Corollary 4.** *Let $p \geq 1$ and $0 \leq m \leq x_i \leq M < \infty$, $p_i > 0$ $(i = 1, .., n)$. Then*

$$(3.4) \qquad 0 \leq \sum_{i=1}^{n} p_i x_i^p - \left(\sum_{i=1}^{n} p_i x_i\right)^p \leq \frac{p}{4}(M-m)\left(M^{p-1} - m^{p-1}\right).$$

*Equality holds iff $x_1 = ... = x_n$.*

The proof follows by (2.17) for the mapping $f(x) = x^p$, $p \geq 1$, $x \geq 0$.

Finally, we have

**Corollary 5.** *Let $p_i, x_i$ be as in Corollary 3. Then*

$$(3.5) \qquad 1 \leq \frac{\prod_{i=1}^{n} x_i^{p_i x_i}}{[A_n(\bar{p}, x)]^{A_n(\bar{p}, x)}} \leq \left(\frac{M}{m}\right)^{\frac{1}{4}(M-m)}.$$

*Equality holds iff $x_1 = ... = x_n$.*

The proof follows by (2.1) for the mapping $f(x) = x \ln x, x > 0$.

## 4. Applications for Shannon's Entropy

Let $X$ be a random variable with the range $R = \{x_1, ..., x_n\}$ and the probability distribution $p_1, ..., p_n$ $(p_i > 0, i = 1, ..., n)$. Define the Shannon entropy mapping

$$H(X) := -\sum_{i=1}^{n} p_i \ln p_i.$$

The following theorem is well known in the literature and concerns the maximum possible value of $H(X)$ in terms of the size of $R$ [15, p. 27].

**Theorem 2.** *Let $X$ be defined as above. Then*

$$(4.1) \qquad 0 \leq H(X) \leq \ln n.$$

*Furthermore, $H(X) = 0$ iff $p_i = 1$ for some $i$ and $H(X) = \ln n$ iff $p_i = \frac{1}{n}$ for all $i \in \{1, ..., n\}$.*

In the recent paper [14], Dragomir and Goh proved the following counterpart result.

**Theorem 3.** *Let $X$ be defined as above. Then*

$$(4.2) \qquad 0 \leq \ln n - H(X) \leq \sum_{1 \leq i < j \leq n} (p_i - p_j)^2.$$

*Equality holds simultaneously in both inequalities iff $p_i = \frac{1}{n}$ for all $i \in \{1, ..., n\}$.*

We give here the following lemma which concerns an analytic inequality for the $\log$-map with applications for entropy.

**Lemma 1.** *Let $0 < m \leq \xi_i \leq M < \infty$, $p_i > 0$ $(i = 1, ..., n)$ with $\sum_{i=1}^{n} p_i = 1$. Then*

$$(4.3) \qquad 0 \leq \ln\left(\sum_{i=1}^{n} p_i \xi_i\right) - \sum_{i=1}^{n} p_i \ln \xi_i \leq \frac{(M-m)^2}{4mM}.$$

The proof follows by Corollary 1, choosing $f(x) = -\ln x$ and $x_i = \xi_i$, $i = 1, ..., n$.

The above lemma provides the following converse inequality for the Shannon entropy mapping (see also [17, p. ]):

**Theorem 4.** *Let $X$ be as above and put $p := \min_{i=\overline{1,n}} p_i$ and $P := \max_{i=\overline{1,n}} p_i$. Then we have*

$$(4.4) \qquad 0 \leq \ln n - H(X) \leq \frac{(P-p)^2}{4pP}.$$

*Proof.* Choose in the above lemma $\xi_i = \frac{1}{p_i} \in \left[\frac{1}{P}, \frac{1}{p}\right]$ and $m = \frac{1}{P}, M = \frac{1}{p}$ to get the desired inequality (4.4). ∎

Another analytic inequality which can be applied for the entropy mapping is embodied in the following lemma.

**Lemma 2.** *Let $0 < m \leq \xi_i \leq M < \infty$, $p_i > 0$ $(i = 1, ..., n)$ with $\sum_{i=1}^{n} p_i = 1$. Then*

$$(4.5) \qquad \begin{aligned} 0 &\leq \sum_{i=1}^{n} p_i \xi_i \ln \xi_i - \sum_{i=1}^{n} p_i \xi_i \ln\left(\sum_{i=1}^{n} p_i \xi_i\right) \\ &\leq \frac{1}{4}(M-m)(\ln M - \ln m) \\ &\leq \frac{1}{4} \cdot \frac{(M-m)^2}{\sqrt{mM}}. \end{aligned}$$

*Proof.* The first inequality follows by Corollary 1, choosing $f(x) = x \ln x$, which is a convex mapping on $(0, \infty)$, and $x_i = \xi_i$, $i = 1, ..., n$.

The second inequality follows by the celebrated inequality between the *geometric mean* $G(a,b) := \sqrt{ab}$ and the *logarithmic mean*

$$L(a,b) := \begin{cases} a & \text{if } b = a \\ \frac{b-a}{\ln b - \ln a} & \text{if } b \neq a \end{cases}, \quad a, b > 0$$

which states that

$$G(a,b) \leq L(a,b), \quad a, b > 0$$

i.e.,

$$\frac{\ln b - \ln a}{b - a} \leq \frac{1}{\sqrt{ab}}.$$

Choosing $b = M$, $a = m$, we obtain

$$\frac{1}{4}(M-m)(\ln M - \ln m) \leq \frac{1}{4}\frac{(M-m)(\ln M - \ln m)}{\sqrt{Mm}}$$

$$= \frac{1}{4\sqrt{Mm}}(M-m)^2.$$

∎

The lemma provides the following converse inequality for the entropy mapping [18].

**Theorem 5.** *Let $X$ be as in Theorem 4. Then we have*

(4.6) $$0 \leq \ln n - H(X) \leq \frac{n}{4}(P-p)(\ln P - \ln p) \leq \frac{n}{4} \cdot \frac{(P-p)^2}{\sqrt{pP}}.$$

*Proof.* Firstly, let us chose $p_i = \frac{1}{n}$ in (4.5) to get

(4.7) $$0 \leq \frac{1}{n}\sum_{i=1}^{n}\xi_i \ln \xi_i - \frac{1}{n}\sum_{i=1}^{n}\xi_i \ln\left(\frac{1}{n}\sum_{i=1}^{n}\xi_i\right)$$

$$\leq \frac{1}{4}(M-m)(\ln M - \ln m)$$

$$\leq \frac{1}{4\sqrt{Mm}}(M-m)^2.$$

Now, if in (4.7) we assume that $\xi_i = p_i \in [p, P]$, then we obtain

$$0 \leq \frac{1}{n}\ln n - \frac{1}{n}H(X) \leq \frac{1}{4}(P-p)(\ln P - \ln p) \leq \frac{1}{4} \cdot \frac{(P-p)^2}{\sqrt{pP}},$$

from where results (4.6). ∎

## 5. Applications for Rényi's entropy

The Rényi entropy of order $\alpha$, $\alpha \in (0,1) \cup (1, \infty)$ is defined by [19]

(5.1) $$H_\alpha(X) := \frac{1}{1-\alpha}\ln\left(\sum_{i=1}^{n}p_i^\alpha\right).$$

Using Jensen's inequality for convex mappings applied for $f(x) = -\ln(x)$, we have

(5.2) $$\ln\left(\sum_{i=1}^{n}p_i x_i\right) \geq \sum_{i=1}^{n}p_i \ln x_i, \ x_i, \ p_i > 0 \ (i=1,...,n), \ \sum_{i=1}^{n}p_i = 1.$$

If we choose $x_i := p_i^{\alpha-1}$ $(i = 1,...,n)$ in (5.2), we obtain

$$\ln\left(\sum_{i=1}^{n}p_i^\alpha\right) \geq (\alpha-1)\sum_{i=1}^{n}p_i \ln p_i,$$

which is equivalent to

(5.3) $$(1-\alpha)[H_\alpha(X) - H(X)] \geq 0.$$

Now, if $\alpha \in (0,1)$, then $H_\alpha(X) \leq H(X)$ and if $\alpha > 1$, then $H_\alpha(X) \geq H(X)$.

Equality holds in (5.3) iff $(p_i)_{i=\overline{1,n}}$ is a uniform distribution and this fact follows by the strict convexity of $-\ln(\cdot)$.

We can now point out a counterpart result for (5.3).

**Theorem 6.** *With the above assumptions, we have*

$$(5.4) \qquad (1-\alpha)\left[H_\alpha(X) - H(X)\right] \leq \frac{1}{4} \frac{\left(P^{\alpha-1} - p^{\alpha-1}\right)^2}{p^{\alpha-1} P^{\alpha-1}},$$

*provided that* $0 < p \leq p_i \leq P < 1$ $(i = 1, ..., n)$.

*Proof.* We use Lemma 1 for $\xi_i := p_i^{\alpha-1}$ and take into account that, for $\alpha \in (0,1)$, we have $P^{\alpha-1} \leq \xi_i \leq p^{\alpha-1}$, $(i = 1, ..., n)$ and, for $\alpha \in (1, \infty)$, we have $p^{\alpha-1} \leq \xi_i \leq P^{\alpha-1}$ $(i = 1, ..., n)$. Choosing in the first case $m = P^{\alpha-1}$, $M = p^{\alpha-1}$ and in the second case $m = p^{\alpha-1}$, $M = P^{\alpha-1}$, we obtain the same upper bound

$$\frac{(M-m)^2}{mM} = \frac{\left(P^{\alpha-1} - p^{\alpha-1}\right)^2}{p^{\alpha-1} P^{\alpha-1}}.$$

∎

Now, let us remark that a particular case of (4.3) for $p_i = \frac{1}{n}$ $(i = 1, ..., n)$ states that

$$(5.5) \qquad 0 \leq \ln\left(\frac{1}{n}\sum_{i=1}^n \xi_i\right) - \frac{1}{n}\sum_{i=1}^n \ln(\xi_i) \leq \frac{(M-m)^2}{4mM},$$

provided that $0 < m \leq \xi_i \leq M < \infty$ $(i = 1, ..., n)$.

This inequality allows us to prove the following inequality for the Rényi entropy:

**Theorem 7.** *With the assumptions of Theorem 6, we have the inequality*

$$(5.6) \qquad \begin{aligned} 0 &\leq (1-\alpha) H_\alpha(X) - \ln n - \alpha \ln G_n(p) \\ &\leq \frac{n}{4} \cdot \frac{(P^\alpha - p^\alpha)^2}{p^\alpha P^\alpha}, \end{aligned}$$

*where* $G_n(p)$ *is the geometric mean of* $p_i$ $(i = 1, ..., n)$, *i.e., we recall that* $G_n(p) = \left(\prod_{i=1}^n p_i\right)^{\frac{1}{n}}$.

*Proof.* We choose in the inequality (5.8) $\xi_i = np_i^\alpha$ $(i = 1, ..., n)$ and observe that $np^\alpha \leq \xi_i \leq nP^\alpha$, so we can have $m = np^\alpha$ and $M = nP^\alpha$ in (5.5), getting the desired inequality. ∎

If we assume that $\alpha \in (0,1)$ and apply Corollary 1 for the convex mapping $f(x) = -x^\alpha$, we deduce the following inequality

$$(5.7) \qquad 0 \leq \left(\sum_{i=1}^n p_i x_i\right)^\alpha - \sum_{i=1}^n p_i x_i^\alpha \leq \frac{\alpha}{4}(M-m)\left(m^{\alpha-1} - M^{\alpha-1}\right),$$

provided that $0 < m \leq x_i \leq M < \infty$, $p_i > 0$ $(i = 1, ..., n)$ and $\sum_{i=1}^n p_i = 1$.

If in (5.7) we put $p_i = \frac{1}{n}$ $(i = 1, ..., n)$, then we obtain

$$(5.8) \qquad 0 \leq \frac{1}{n^\alpha}\left(\sum_{i=1}^n x_i\right)^\alpha - \frac{1}{n}\sum_{i=1}^n x_i^\alpha \leq \frac{\alpha}{4}(M-m)\left(m^{\alpha-1} - M^{\alpha-1}\right).$$

The following inequality for the $\alpha$-Rényi entropy holds.

**Theorem 8.** *If $0 < p \leq p_i \leq P < 1$ $(i = 1, ..., n)$ and $\alpha \in (0, 1)$, then we have the inequality*

$$(5.9) \qquad 0 \leq n^{1-\alpha} - \exp\left[(1 - \alpha) H_\alpha(X)\right] \leq \frac{\alpha}{4} n (P - p) \left(p^{\alpha-1} - P^{\alpha-1}\right).$$

*Proof.* Choosing $x_i = p_i$ $(i = 1, ..., n)$ in (5.8), we deduce

$$(5.10) \qquad 0 \leq \frac{1}{n^\alpha} - \frac{1}{n} \sum_{i=1}^n p_i^\alpha \leq \frac{\alpha}{4} (P - p) \left(p^{\alpha-1} - P^{\alpha-1}\right).$$

Taking into account that $\sum_{i=1}^n p_i^\alpha = \exp\left[(1 - \alpha) H_\alpha(X)\right]$, then from (5.10) we deduce the desired inequality (5.9). ∎

Now, if we define by $E(X) := \sum_{i=1}^n p_i^2$, the informational energy of the random variable X, then we may also point out the following inequality.

**Theorem 9.** *If $p_i$ and $\alpha$ are as in Theorem 8, then we have the inequality*

$$(5.11) \qquad 0 \leq E^\alpha(X) - \exp\left[-\alpha H_{\alpha+1}(X)\right] \leq \frac{\alpha}{4} (P - p) \left(p^{\alpha-1} - P^{\alpha-1}\right).$$

The proof follows by the inequality (5.7), choosing $x_i = p_i$, $i = 1, ..., n$, and we omit the details.

Now assume that $\alpha \in (1, \infty)$. Then, applying Corollary 1 for the convex mapping $f(x) = x^\alpha$, we deduce the following inequality

$$(5.12) \qquad 0 \leq \sum_{i=1}^n p_i x_i^\alpha - \left(\sum_{i=1}^n p_i x_i\right)^\alpha \leq \frac{\alpha}{4} (M - m) \left(m^{\alpha-1} - M^{\alpha-1}\right),$$

provided that $0 < m \leq x_i \leq M < \infty$, $p_i > 0$ $(i = 1, ..., n)$ and $\sum_{i=1}^n p_i = 1$.

Finally, using a similar argument to the one above, we can state the following theorem.

**Theorem 10.** *Let $\alpha \in (1, \infty)$ and $0 < p \leq p_i \leq P < 1$ $(i = 1, ..., n)$. Then we have the inequalities*

$$0 \leq \exp\left[-\alpha H_\alpha(X)\right] - E^\alpha(X) \leq \frac{\alpha}{4} (P - p) \left(P^{\alpha-1} - p^{\alpha-1}\right)$$

*and*

$$0 \leq \exp\left[(1 - \alpha) H_\alpha(X)\right] - n^{1-\alpha} \leq \frac{\alpha}{4} n (P - p) \left(P^{\alpha-1} - p^{\alpha-1}\right).$$

School of Communications and Informatics, Victoria University of Technology, PO Box 14428, Melbourne City MC 8001, Australia

*E-mail address*: `sever@matilda.vu.edu.au`